\documentclass[12]{amsart}
\usepackage{amsfonts}
\usepackage{mathrsfs}

\usepackage{amsmath}
\usepackage{amssymb}
\usepackage{graphicx}
\usepackage{geometry}
\usepackage{enumerate}
\newtheorem{thm}{Theorem}[section]
\newtheorem{cor}[thm]{Corollary}
\newtheorem{lem}[thm]{Lemma}

\newtheorem{prop}[thm]{Proposition}
\newtheorem{exa}[thm]{Example}

\theoremstyle{definition}
\newtheorem{defn}[thm]{Definition}
\newtheorem{rem}[thm]{Remark.}

\numberwithin{equation}{section}

\begin{document}
\vbox{\vskip 1cm}
\title{On a question of Mohamed and M\"uller} \maketitle
\begin{center}
Yongduo Wang\ \ \ \ Dejun Wu \ \ \\
\vskip 4mm Department of Applied Mathematics, Lanzhou University of Technology\\
Lanzhou 730050, P. R. China\\
 E-mail: ydwang@lut.cn\\
E-mail: wudj2007@gmail.com \\

\end{center}

\bigskip

\noindent {\bf Abstract:} A module $M$ is called H-supplemented if
for every submodule $A$ of $M$ there is a direct summand $A'$ of $M$
such that $A+X=M$ holds if and only if $A'+X=M$ for any submodule
$X$ of $M$. (Equivalently, for each $X\leq M$, there exists a direct
summand $D$ of $M$ such that $(X+D)/D\ll M/D$ and $(X+D)/X\ll M/X$.)
Direct summands and sums of H-supplemented modules are studied and a
question posed by Mohamed and M\"uller in 1990 is answered in the
negative.


\vskip 2mm \noindent {\bf MR(2000)Subject Classification:} 16D10;
16D99

\bigskip

\baselineskip=20pt

\section{Introduction}

\noindent Lifting modules and extending modules, as dual concepts,
play important roles in rings and categories of modules. Lifting
modules and their generalizations have been studied extensively by
many authors recently. As a proper generalization of lifting
modules, the notion of H-supplemented modules was introduced in
Mohamed and M\"uller (1990). Let $M$ be a module. $M$ is called
\emph{H-supplemented} if for every submodule $A$ there is a direct
summand $A'$ of $M$ such that $A+X=M$ holds if and only if $A'+X=M$
for any submodule $X$ of $M$. Equivalently, for each $X\leq M$,
there exists a direct summand $D$ of $M$ such that $(X+D)/D\ll M/D$
and $(X+D)/X\ll M/X$. It is not known whether the class of
H-supplemented modules is closed under direct summands. $M$ is
called \emph{completely H-supplemented} if every direct summand of
$M$ is H-supplemented. Direct summands and sums of H-supplemented
modules are studied and a question (Is every H-supplemented module
amply supplemented?) posed by Mohamed and M\"uller in 1990 is
answered in the negative in this paper.

Throughout this paper all rings are associative with unity and all
modules will be unital right $R$-modules. We use $N\leq M$ ($N\ll
M$) to indicate that $N$ is a submodule (a small submodule) of $M$.
Other terminology and notation can be found in Anderson et al.
(1992), Clark et al. (2006) and Mohamed and M\"uller (1990).

\section{Main Results}

\indent Let $M$ be a module. Recall the following conditions:

 $(D_2)$: If $N$ is a submodule of $M$ such that $M/N$ is isomorphic to a direct summand
of $M$, then $N$ is a direct summand of $M$.

 $(D_3)$: For every direct summands $K,L$ of $M$ with $M=K+L, K\cap L$ is a direct
summand of $M$.

\begin{thm}
Let $M$ be an H-supplemented module. If $M$ has $(D_3)$, then $M$ is
completely H-supplemented.
\end{thm}

\begin{proof}
Let $N$ be a direct summand of $M$. Then $M=N\oplus N', N'\leq M$.
Let $Y$ be a submodule of $N$. Since $M$ is H-supplemented, there
exists a submodule $X$ of $M$ and a direct summand $K$ of $M$ such
that $X/(Y\oplus N')\ll M/(Y\oplus N')$ and $X/K\ll M/K$. Then
$M=X+N$ and so $M=K+N$. Now $K\cap N$ is a direct summand of $M$ as
$M$ has $(D_3)$. Hence $K\cap N$ is a direct summand of $N$. Note
that $(X\cap N)/(K\cap N)\ll N/(K\cap N)$ and $(X\cap N)/Y\ll N/Y$
by Ganesan and Vanaja (2002, Lemma 2.6) and Keskin (2000, Lemma
1.1). Thus $N$ is H-supplemented. Therefore, $M$ is completely
H-supplemented.
\end{proof}
\begin{exa} Let $R$ be a discrete valuation ring with field of
fractions $K$. Let $P$ be a unique maximal ideal of $R$. Then the
$R$-module $M=(K/R)\oplus (R/P)$ is an H-supplemented module with
$(D_{3})$. Hence $M$ is completely H-supplemented by Theorem 2.1.
\end{exa}

\begin{exa} Let $R$ be a discrete valuation ring with field of fractions $K$.
Assume $K_R$ is not quasi-projective (and hence $R$ is a Dedekind
domain). By Keskin (2000, Example 3.7, 2002, Example 3.8),
$R$-module $M=K\oplus K$ is a $\oplus$-supplemented module with
$(D_{3})$, but not amply supplemented. Note that over Dedekind
domains, H-supplemented modules and $\oplus$-supplemented modules
are the same by Mohamed and M\"uller (1990, Proposition A.7 and
A.8). Hence $M$ is H-supplemented. Therefore, it is completely
H-supplemented by Theorem 2.1.
\end{exa}

\begin{rem} Example 2.3 shows that there is an H-supplemented module
which is not amply supplemented. Thus we answer an open question (Is
every H-supplemented module amply supplemented?) posed by Mohamed
and M\"uller (1990, P106).
\end{rem}

The following examples show that the condition $(D_3)$ in Theorem
2.1 is not necessary.
\begin{exa} It is well known that $\mathbb{Z}/2\mathbb{Z}\oplus
\mathbb{Z}/8\mathbb{Z}$ is completely H-supplemented, but not
$(D_3)$.
\end{exa}
\begin{exa} It is well known that $\mathbb{Z}/p\mathbb{Z}\oplus
\mathbb{Z}/p^{2}\mathbb{Z}$ ($p$ is prime) is a lifting module, and
hence it is a completely H-supplemented module. However, it does not
satisfy the condition $(D_3)$.
\end{exa}

\begin{cor}
Let $M$ be an H-supplemented module with $(D_2)$ or $(D_3)$. Then so
is every direct summand of $M$.
\end{cor}

\begin{proof}
Clear.
\end{proof}

Recall that a module $M$ is a \emph{UC} module if every submodule of
$M$ has a unique closure in $M$.

\begin{cor}
Let $M$ be an extending UC module. Then $M$ is H-supplemented if and
only if $M$ is completely H-supplemented.
\end{cor}

\begin{proof}
Let $K$ and $L$ be direct summands of $M$ with $K+L=M$. Since $M$ is
a UC module, $K\cap L$ is a closed submodule of $M$. Then $K\cap L$
is a direct summand of $M$. Now the result follows by Theorem 2.1.
\end{proof}

Recall that a module $M$ is called a \emph{polyform} module if every
essential submodule of $M$ is a rational submodule of $M$.

\begin{cor}
Let $M$ be a polyform extending module. Then $M$ is H-supplemented
if and only if $M$ is completely H-supplemented.
\end{cor}

\begin{proof}
It is well known that a polyform module is a UC module. Now the
result follows by Corollary 2.8.
\end{proof}

\begin{cor}
Let $M$ be a quasi-injective polyform module. Then for each
$n\in\mathbb{N}$,  $M^n$ is H-supplemented if and only if $M^n$ is
completely H-supplemented.
\end{cor}

\begin{proof}
Let $M^n$ be an H-supplemented module. Since $M$ is quasi-injective
polyform module, $M^n$ is quasi-injective polyform module. Now the
result follows by Corollary 2.9.
\end{proof}

\begin{defn} (See Wu and Wang, Definition 3.1)
Let $M_1$ and $M_2$ be modules such that $M=M_1\oplus M_2$. We say
$M_1$ is $M_2$-sjective if for every $A\leq M$ such that $M=A+M_2$,
there exists $K\leq M$ such that $M=K\oplus M_2$ and $(A+K)/A\ll
M/A$. $M_1$ and $M_2$ are called relatively sjective if $M_1$ is
$M_2$-sjective and $M_2$ is $M_1$-sjective.
\end{defn}

Clearly, if $M_1$ is $M_2$-projective, then $M_1$ is $M_2$-sjective.
The converse need not be true in general by Example 2.13.

\begin{lem} Let $M=M_1\oplus M_2$ be an H-supplemented module with $(D_3)$. Then
$M_1$ and $M_2$ are relatively sjective.
\end{lem}

\begin{proof}
See Wu and Wang (Proposition 3.2).
\end{proof}

\begin{exa} Let $R$ be a discrete valuation ring with field of
fractions $K$. Let $P$ be a unique maximal ideal of $R$. Then the
$R$-module $M=(K/R)\oplus (R/P)$ is an H-supplemented module with
$(D_{3})$. Hence $K/R$ and $R/P$ are relatively sjective by Lemma
2.12. However, $K/R$ and $R/P$ are not relatively projective (If
$K/R$ and $R/P$ are relatively projective. Since $M$ is amply
supplemented and $K/R$ and $R/P$ are lifting, $M$ is lifting by
Keskin (2000, Corollary 2.9). But $M$ is not lifting. This is a
contradiction.) So this example shows that Definition 2.11 is a
proper generalization of projective modules (Also see Example 2.3).
\end{exa}

\begin{lem}
Let $M=M_1\oplus M_2$ be a weakly supplemented module. If $M_1$ is
$M_2$-sjective (or $M_2$ is $M_1$-sjective) and $M_1$ and $M_2$ are
H-supplemented, then $M$ is H-supplemented.
\end{lem}

\begin{proof}
See Wu and Wang (Theorem 3.4).
\end{proof}

Recall that a module $M$ is called a \emph{FI-lifting module}, if
for any fully invariant submodule $N$ of $M$, there exists a direct
summand $K$ of $M$ such that $K\leq N$ and $N/K\ll M/K$.

\begin{thm}
Let $M$ be a weakly supplemented module with $(D_3)$. Then the
following statements are equivalent.
\begin{enumerate}[(i)]
\item $M$ is completely H-supplemented.
\item $M$ is H-supplemented.
\item $M=K\oplus K'$, where $K$ and $K'$ are H-supplemented, Rad$(M)/K\ll
M/K$, Rad$(K')\ll K'$ and $K$ and $K'$ are relatively sjective.
\item $M=K\oplus K'$, where $K$ and $K'$ are H-supplemented, Soc$(M)/K\ll
M/K$, Soc$(K')\ll K'$ and $K$ and $K'$ are relatively sjective.
\item $M=K\oplus K'$, where $K$ and $K'$ are H-supplemented, $Z(M)/K\ll
M/K$, $Z(K')\ll K'$ and $K$ and $K'$ are relatively sjective.
\end{enumerate}
\end{thm}

\begin{proof}
$(i)\Leftrightarrow(ii)$ By Theorem 2.1.

$(ii)\Rightarrow(iii)$ Let $M$ be an H-supplemented module. Then it
is easy to see that $M$ is FI-lifting. There exists a direct summand
$K$ of $M$ such that Rad$(M)/K\ll M/K$. Write $M=K\oplus K', K'\leq
M$. Clearly, Rad$(M)=$Rad$(K)\oplus$Rad$(K')$ and so
Rad$(K')=K'\cap$Rad$(M)$. Hence Rad$(K')\ll K'$. By Lemma 2.12, $K$
and $K'$ are relatively ejective.

$(iii)\Rightarrow(ii)$ By Lemma 2.14.

$(ii)\Leftrightarrow(iv)$ and $(ii)\Leftrightarrow(v)$ are similar
to $(ii)\Leftrightarrow(iii)$.
\end{proof}

\begin{prop}
Let $M=M_1\oplus M_2$. If for every submodule $N$ of $M_1$ there
exists a direct summand $K$ of $M$ such that $M_2\leq K$,
$(N+K)/K\ll M/K$ and $(N+K)/N\ll M/N$, then $M_1$ is H-supplemented.
\end{prop}

\begin{proof}
Let $L$ be a submodule of $M_1$. By hypothesis, there exists a
direct summand $K$ of $M$ such that $M_2\leq K$, $(L+K)/K\ll M/K$
and $(L+K)/L\ll M/L$. Now $K=(K\cap M_1)\oplus M_2$. Hence $K\cap
M_1$ is a direct summand of $M_1$. It is easy to check that
$(L+(K\cap M_1))/L\ll M_1/L$ and $(L+(K\cap M_1))/(K\cap M_1)\ll
M/(K\cap M_1)$. Therefore, $M_1$ is H-supplemented.
\end{proof}

\begin{thm}
Let $M=M_1\oplus M_2$, where $M_1$ is a fully invariant submodule of
$M$. If $M$ is H-supplemented, then $M_1$ and $M_2$ are
H-supplemented.
\end{thm}

\begin{proof}
By Kosan and Keskin (2007, Corollary 2.4), $M_2$ is H-supplemented.
Next we show that $M_1$ is H-supplemented. Let $K$ be a submodule of
$M_1$. Since $M$ is H-supplemented, there exists a direct summand
$D$ of $M$ such that $(K+D)/K\ll M/K$ and $(K+D)/D\ll M/D$. Write
$M=D\oplus D', D'\leq M$. Then $M=K+D'$. Since $M_1$ is a fully
invariant submodule of $M$, $M_1=(M_1\cap D)\oplus(M_1\cap D')$.
Hence $M=M_1+D'=(M_1\cap D)\oplus D'$. Thus $D=D\cap M_1$ and so
$D\leq M_1$. Now $(K+D)/K\ll M_1/K$ and $(K+D)/D\ll M_1/D$.
Therefore, $M_1$ is H-supplemented.
\end{proof}

Let $R$ be any ring and $M$ be an $R$-module. The module $M$ is
called small if $M\ll E(M)$, where $E(M)$ is the injective hull of
$M$. Talebi and Vanaja (2002) defined $\bar{Z}(M)=\cap\{\text{Ker}g
| g : M\rightarrow N \ \text{and}\  N \ \text{is small}\}$. They
call $M$ cosingular if $\bar{Z}(M)=0$ and noncosingular if
$\bar{Z}(M)=M$. Let $M$ be a module. Talebi and Vanaja defined
$\bar{Z}^0(M)= M, \bar{Z}^1(M) =\bar{Z}(M)$ and defined inductively
$\bar{Z}^{\alpha}(M)$ for any ordinal $\alpha$. Thus, if $\alpha$ is
not a limit ordinal they set
$\bar{Z}^{\alpha}(M)=\bar{Z}(\bar{Z}^{\alpha-1}(M))$, while if
$\bar{Z}^{\alpha}(M)$ is a limit ordinal they set
$\bar{Z}^{\alpha}(M)=\cap_{\beta<\alpha}\bar{Z}^{\beta}(M)$. This
gives the descending sequence
$M=\bar{Z}^0(M)\supseteq\bar{Z}^1(M)\supseteq\bar{Z}^2(M)...$ of
submodules of M (see Talebi and Vanaja (2002)). Let $M$ be a module.
It is proved that in Talebi and Vanaja (2002, Theorem 4.1) that $M$
is lifting if and only if $M=\bar{Z}^2(M)\oplus M', M'\leq M$ such
that $M'$ and $\bar{Z}^2(M)$ are lifting, $M'$ is
$\bar{Z}^2(M)$-projective and $M$ is amply supplemented. Here we
have the following fact.

\begin{thm}
Let $M$ be an amply supplemented module. Then $M$ is H-supplemented
if and only if $M=\bar{Z}^2(M)\oplus M'$, where $\bar{Z}^2(M)$ and
$M'$ are H-supplemented.
\end{thm}

\begin{proof}
Let $M$ be an H-supplemented module. Note that $\bar{Z}^2(M)$ is a
fully invariant coclosed submodule of $M$. Since $M$ is FI-lifting,
$M=\bar{Z}^2(M)\oplus M'$, where $\bar{Z}^2(M)$ and $M'$ are
H-supplemented by Theorem 2.17. Conversely, let $M$ be an amply
supplemented module. Then $\bar{Z}^2(M)$ is $M'$-projective by the
proof of Talebi and Vanaja (2002, Theorem 4.1). Therefore, $M$ is
H-supplemented by Lemma 2.14.
\end{proof}

It is well known that a direct sum of H-supplemented modules need
not be H-supplemented. Here we will study an infinite direct sum of
H-supplemented modules. Let $M$ be an $R$-module such that
$M=\oplus_{i\in I}M_i$ is the direct sum of H-supplemented modules
$M_i(i\in I)$, for some given index set $I$. Now we consider when
$M$ itself is an H-supplemented module. Let $M=\oplus_{i\in I}M_i$.
For each $i\in I$, $M_{-i}$ will denote $\oplus_{j\in
I\backslash\{i\}}M_j$. For any set $I$, $|I|$ will denote its
cardinality.

\begin{thm}
Let $M=\oplus_{i\in I}M_i$ be the direct sum of modules $M_i(i\in
I)$, for some index set $I$ with $|I|\geq 2$. If $M$ is a weakly
supplemented module with $(D_3)$, then the following statements are
equivalent.
\begin{enumerate}[(i)]
\item $M$ is H-supplemented.
\item There exists $i\in I$ such that for
every submodule $K$ of $M$ with $M=K+M_i$ or $M=K+M_{-i}$ there
exists a direct summand $N$ of $M$ such that $(K+N)/K\ll M/K$ and
$(K+N)/N\ll M/N$.
\item There exists $i\in
I$ such that for every submodule $K$ of $M$ with $(K+M_i)/K\ll M/K$
or $(K+M_{-i})/K\ll M/K$ or $M=K+M_i=K+M_{-i}$ there exists a direct
summand $N$ of $M$ such that $(K+N)/K\ll M/K$ and $(K+N)/N\ll M/N$.
\end{enumerate}
\end{thm}

\begin{proof}
$(i)\Rightarrow(ii)$ is clear.

$(ii)\Rightarrow(i)$ Let $N\leq M$. If $M=N+M_i$, then there is
nothing to prove. Now assume that $M\neq N+M_i$. Consider the
submodule $(N+M_i)/N$ of $M/N$. Since $M$ is weakly supplemented,
there exists a submodule $K/N$ of $M/N$ such that
$(N+M_i)/N+K/N=M/N$ and $((N+M_i)\cap K)/N\ll M/N$. Then $M=K+M_i$.
By hypothesis, there exists a direct summand $L$ of $M$ such that
$(K+L)/K\ll M/K$. Clearly, $M=(N+M_i)+M_{-i}$. By (ii), there exists
a submodule $X_1$ of $M$ and a direct summand $D_1$ of $M$ such that
$X_1/(N+M_i)\ll M/(N+M_i)$ and $X_1/D_1\ll M/D_1$. Obviously,
$M=(K+L)+M_i$. There exists a submodule $X_2$ of $M$ and a direct
summand $D_2$ of $M$ such that $X_2/(K+L)\ll M/(K+L)$ and
$X_2/D_2\ll M/D_2$. Then $X_2/K\ll M/K$ and so $(X_1\cap
X_2)/((N+M_i)\cap K)\ll M/((N+M_i)\cap K)$. Therefore $(X_1\cap
X_2)/N\ll M/N$. Note that $M=X_1+X_2$. Hence $M=D_1+D_2$. Thus
$(X_1\cap X_2)/(D_1\cap D_2)\ll M/(D_1\cap D_2)$. Since $M$ has
$(D_3)$, $D_1\cap D_2$ is a direct summand of $M$. Therefore, $M$ is
H-supplemented.

$(ii)\Rightarrow(iii)$ is obvious.

$(iii)\Rightarrow(ii)$ Let $N\leq M$ with $M=N+M_i$, the case
$M=N+M_{-i}$ being analogous. If $M=N+M_{-i}$, then there is nothing
to prove. Now assume that $M\neq N+M_{-i}$. Consider the submodule
$(N+M_{-i})/N$ of $M/N$. Since $M$ is weakly supplemented, there
exists a submodule $K/N$ of $M/N$ such that $(N+M_{-i})/N+K/N=M/N$
and $((N+M_{-i})\cap K)/N\ll M/N$. Then $M=K+M_i=K+M_{-i}$. By
hypothesis, there exists a direct summand $L$ of $M$ such that
$(K+L)/K\ll M/K$. Clearly, $((N+M_{-i})+M_{-i})/(N+M_{-i})\ll
M/(N+M_{-i})$. By (iii), there exists a submodule $X_1$ of $M$ and a
direct summand $D_1$ of $M$ such that $X_1/(N+M_{-i})\ll
M/(N+M_{-i})$ and $X_1/D_1\ll M/D_1$. Obviously,
$M=(K+L)+M_i=(K+L)+M_{-i}$. There exists a submodule $X_2$ of $M$
and a direct summand $D_2$ of $M$ such that $X_2/(K+L)\ll M/(K+L)$
and $X_2/D_2\ll M/D_2$. Then $(X_1\cap X_2)/N\ll M/N$. Note that
$M=D_1+D_2$. Thus $(X_1\cap X_2)/(D_1\cap D_2)\ll M/(D_1\cap D_2)$.
Since $M$ has $(D_3)$, $D_1\cap D_2$ is a direct summand of $M$. Now
the proof is completed.
\end{proof}

\end{document}